\documentclass[review]{elsarticle}

\usepackage{lineno,hyperref}
\modulolinenumbers[5]

\usepackage{amsmath}
\usepackage{amssymb}
\usepackage{graphicx} % figures
\usepackage{algorithm}
\usepackage{algpseudocode}
\renewcommand{\vec}[1]{\boldsymbol{#1}} % for BOLD vectors.

\journal{arXiv.org}

\begin{document}

\begin{frontmatter}

\title{Global finite element matrix construction based on a CPU-GPU implementation}

%% Group authors per affiliation:
\author[1address]{Francisco Javier Ram\'{i}rez-Gil}
\ead{fjramireg@unal.edu.co}

\author[2address]{Marcos de Sales Guerra Tsuzuki}
\ead{mtsuzuki@usp.br}

\author[1address]{Wilfredo Montealegre-Rubio \corref{correspondingauthor}}
\ead{wmontealegrer@unal.edu.co}
\cortext[correspondingauthor]{Corresponding author}

\address[1address]{Department of Mechanical Engineering, Faculty of Mines, Universidad Nacional de Colombia, Medell\'{i}n, 050034, Colombia}
\address[2address]{Computational Geometry Laboratory, Department of Mechatronics and Mechanical Systems Engineering, Escola Polit\'{e}cnica da Universidade de S\~{a}o Paulo,  Av. Prof. Mello Moraes, 2231 S\~{a}o Paulo - SP, Brazil}

\begin{abstract}
The finite element method (FEM) has several computational steps to numerically solve a particular problem, to which many efforts have been directed to accelerate the solution stage of the linear system of equations. However, the finite element matrix construction, which is also time-consuming for unstructured meshes, has been less investigated. The generation of the global finite element matrix is performed in two steps, computing the local matrices by numerical integration and assembling them into a global system, which has traditionally been done in serial computing. This work presents a fast technique to construct the global finite element matrix that arises by solving the Poisson's equation in a three-dimensional domain. The proposed methodology consists in computing the numerical integration, due to its intrinsic parallel opportunities, in the graphics processing unit (GPU) and computing the matrix assembly, due to its intrinsic serial operations, in the central processing unit (CPU). In the numerical integration, only the lower triangular part of each local stiffness matrix is computed thanks to its symmetry, which saves GPU memory and computing time. As a result of symmetry, the global sparse matrix also contains non-zero elements only in its lower triangular part, which reduces the assembly operations and memory usage. This methodology allows generating the global sparse matrix from any unstructured finite element mesh size on GPUs with little memory capacity, only limited by the CPU memory. 
\end{abstract}


\end{frontmatter}

%\linenumbers

%%%% Start %%%%%%
\section{Introduction}
\label{sec1} 
Numerous physical phenomena in a stationary condition such as electrical and magnetic potential, heat conduction, fluid flow and elastic problems in a static condition can be described by elliptic partial differential equations (EPDEs). The EPDEs does not involve a time variable, and so describes the steady state response. A linear EPDE has the general form given in Eq. ({\ref{Eq.1}}), where $ a,b,c $ are coefficient functions, $ f $ is a source (excitation) function and $ u $ is the unknown variable. All of these functions can vary spatially $ (x,y,z) $. 

\begin{equation}
\nabla(c \cdot \nabla u) + b \cdot \nabla u + au = f
\label{Eq.1} 
\end{equation}

EPDEs can be solved exactly by mathematical procedures such as Fourier series {\cite{Kreyszig2011}}. However, the classical solution does not frequently exist; for those problems which allow the use of these analytical methods, many simplifications are made {\cite{Braess2007}}. Consequently, several numerical methods have been developed such as the finite element method (FEM) and finite difference method to efficiently solve EPDEs. 

The FEM has several advantages over other methods. The main advantage is that it is particularly useful for problems with complicated geometry using unstructured meshes {\cite{Braess2007}}. One way to get a suitable framework for solving EPDEs problems, by using FEM, is formulating them as variational problems, also called weak solution.

\begin{figure}[hbtp]
\centering
\includegraphics[width=1\textwidth]{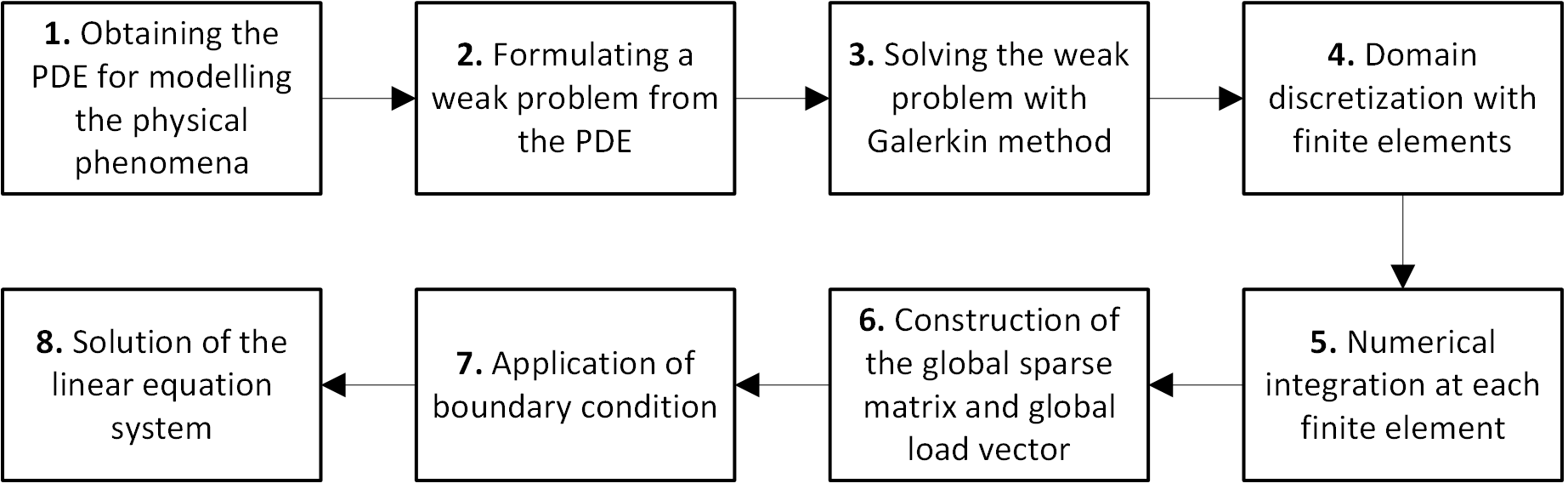}
\caption{FEM steps to solve EPDEs.}
\label{FEM_flow}
\end{figure}

The variational formulation of an EPDE is a mathematical treatment for converting the strong formulation into a weak formulation, which permits the approximation in elements or subdomains, and the EPDE variational form is solved using the so-called Galerking method {\cite{Zienkiewicz2005}}. Fig. {\ref{FEM_flow}} summarizes the steps required to numerically solve an EPDE by using FEM. Some projects are publicly available and all of these steps are solved automatically in the computer {\cite{Logg2012}}; however, the most common commercial FEM packages are focused on steps 4-8, which are: 

\begin{enumerate}[(i)]
    \item Domain discretization with finite elements (FEs);
    \item Numerical integration of all FEs. In this step, a local matrix $\vec{k}_e$ and a local load vector $\vec{f}_e$ are computed for each finite element;
    \item Construction of the global sparse matrix $\vec{K}$ from local matrices $\vec{k}_e$ and global load vector $\vec{F}$ from local load vectors $\vec{f}_e$;
    \item Application of boundary conditions;
    \item Solution of the linear equation system formed previously, $\vec{K} \vec{U}=\vec{F}$, where $\vec{U}$ represents the vector of the nodal solution. 
\end{enumerate}

Many efforts have been made to accelerate the linear equation solver by using direct {\cite{Davis2006}} and iterative methods {\cite{Saad2003}}; however, the FE sparse matrix construction (steps 5 and 6 of Fig. {\ref{FEM_flow}}), which are also time-consuming for dense unstructured meshes, have been less investigated {\cite{Dziekonski2012a}}. Hence, in this work the fast construction of the global sparse matrix is developed, which arises in the numerical solution of EPDEs by FEM. 

One way to efficiently construct the global FE matrix is by using parallel computing. Parallelism is a trend in current computing and, for this reason, microprocessor manufacturers focus on adding cores rather than on increasing single-thread performance {\cite{Owens2008}}. Although there are several parallel architectures, the most widespread many-core architecture are the GPUs (graphics processing units), and the most common multi-core architecture are the CPUs (central processing units) {\cite{Brodtkorb2010}}. However, due to the massive parallelism offered by GPU compared to CPU, the trend in scientific computing is to use GPU for accelerating intensive and intrinsic parallel computations {\cite{Kirk&Hwu2012}}. Additionally, the GPU is preferred over the other architectures to be used as a co-processor to accelerate part or the whole code because its market is well established in the computer games industry. This fact allows having a graphic processor in most personal computers, which makes GPUs an economic and attractive option for application developers {\cite{Kirk&Hwu2012}}.

Although some researches are focused on completely porting the FEM subroutines to be executed in GPU {\cite{Komatitsch2009, Fu2014}}, other works are focused on porting parts of the FEM code to the GPU, such as numerical integration {\cite{Knepley2013, Maciol2010, Kruzel2013, Banas2014}}, matrix assembly 	 {\cite{Cecka2011, Markall2013}}, and the solution of large linear systems of equations  {\cite{JKW2005, Stock2010, Li&Saad2012}}. Even though some previous works accelerated the construction of sparse FE matrix {\cite{Dziekonski2012a, Dziekonski2013, Dziekonski2012b}}, these approaches require high GPU memory, which means that the size of the problem that may be analyzed is limited or an expensive GPU with high memory capacity is needed. However, in most personal computers equipped with a graphics processor, the GPU has low-memory capacity.

Accordingly, the proposed approach is based on the assumption that a single GPU has significantly less memory than the CPU; therefore, a CPU-GPU implementation is developed. The central idea is computing the numerical integration in the GPU, which is an intensive and intrinsic parallel subroutine, and computing the assembly in the CPU, which is a predominantly serial subroutine with high memory consumption. Hence, based on this methodology, we are able to construct any large sparse matrix arising in FEM, only limited by the CPU memory.

The rest of the paper is organized as follows. In section {\ref{sec2}}, the FEM formulation of the Poisson's equation used as an example for checking the proposed implementation is given. Section {\ref{sec3}} presents the CPU-GPU implementation of the FE matrix construction. The results are discussed in section {\ref{sec4}} and finally, in section {\ref{sec5}}, the conclusions are drawn.

\section{Elliptic PDE modelling by FEM}
\label{sec2} 
The computational benefits of our proposed implementation in the construction of large sparse FE matrices is shown by solving the Poisson's equation in a three-dimensional (3D) domain. The Poisson's equation is an elliptic-type PDE. This type of equation is commonly used to solve scalar problems such as electrical and magnetic potential, heat transfer and other problems with no external flows {\cite{Zienkiewicz2005}}. In addition, the Poisson's equation is used for modelling vector field problems such as structural problems without external applied forces. All of these problems should be in static or stationary condition. In this section, the FEM steps are presented, from the strong formulation step to the construction of the global sparse matrix step (steps 1 to 6 of Fig. {\ref{FEM_flow}}). 

The strong formulation of the Poisson equation is given by {\cite{Zienkiewicz2005}}:

\begin{equation} 
\begin{array}{c}
\nabla \cdot (c \nabla \phi(\vec{r})) = 0 \ \ \textup{in} \ \ \Omega \\
\phi = \hat{\phi} \ \ \textup{on} \ \ \Gamma_\phi \\
q_n = 0 \ \ \textup{on} \ \ \Gamma_q
\end{array}
\label{Eq.2} 
\end{equation}
where $ \phi $ represents a general scalar variable on the domain $ \Omega $, i.e. can be the electrical potential or temperature, depending on the specific applications. The position vector is $ \vec{r} $ and the term $ c $ can represent the electrical or thermal conductivity for an isotropic material. Term $ \hat{\phi} $ is a prescribed value of scalar variable $ \phi $ at boundary $ \Gamma_\phi $ (Dirichlet condition) and the normal component of the flow (or flux) $ q_n $ is stated for simplicity as zero at boundary $ \Gamma_q $ (Neumann condition) for the purposes of computational evaluation only. After the integration by parts of Eq. ({\ref{Eq.2}}), the weak form is found as:

\begin{equation}
\displaystyle \int_\Omega (\nabla \nu)^{T} (c \nabla \phi) d\Omega = 0
\label{Eq.3} 
\end{equation}
where $ \nu $ is an arbitrary function {\cite{Zienkiewicz2005}}. Based on the weak form of the Poisson equation, the FE formulation can be found. Accordingly, the domain is discretized into FEs; next, in each FE, the scalar variable is approximated as: 

\begin{equation}
\phi \approx \hat{\phi} = \sum_{a=1}^{n} N_a \tilde{\phi}_a = \vec{N} \tilde{\phi}
\label{Eq.4}
\end{equation}
where $ n $ is the number of nodes in the element and, according to the Galerkin principle, $ \nu = N_a $ are the shape functions. Substituting Eq. ({\ref{Eq.4}}) into Eq. ({\ref{Eq.3}}) and evaluating the integrals for all elements, a system of equation is obtained as:

\begin{equation}
\vec{K} \vec{\Phi} = \vec{F}; \ \ \ \vec{F} = \vec{0}
\label{Eq.5}
\end{equation}
where $ \vec{K} $ is the so-called global sparse matrix or \emph{stiffness} matrix. Term $ \vec{ \Phi} $ is the nodal solution vector and global load vector $ \vec{F} $ is modified by the boundary conditions, ceasing to be a vector of zeroes. The $ \vec{K} $ matrix is obtained as a \emph{sum} of contributions of the matrices obtained in numerical integration for each FE. This step is known as assembly {\cite{Zienkiewicz2005}}: 

\begin{equation}
\vec{K} = \sum_{e=1}^{nel} \vec{k}_e
\label{Eq.6}
\end{equation}
where $ nel $ is the number of FEs and $ \vec{k}_e $ is the local matrix of element $ e $ (local \emph{stiffness} matrix), which is obtained as:  

\begin{equation}
\vec{k}_e = \int_{\Omega^e} c \vec{B}^T \vec{B} d \Omega, \ \ \textup{with} \ \ \vec{B}=\sum_{a=1}^n \nabla N_a
\label{Eq.7}
\end{equation}

In our approach we consider that coefficient $ c $, as a part of material properties array ($ \vec{MP} $), can vary from one element to the other, but it remains constant for a specific element, which means it can be outside the integral in Eq. ({\ref{Eq.7}}). The global sparse matrix construction is made in Eq. ({\ref{Eq.6}}) and Eq. ({\ref{Eq.7}}). However, for computational efficiency, Eq. ({\ref{Eq.7}}) is solved numerically instead of by exact integration. There are several numerical procedures to solve an integral, but the most common in FEM is the Gauss-Legendre, or simply Gauss quadrature {\cite{Zienkiewicz2005}}. For using the Gauss quadrature, a transformation by using the reference element technique must be applied, in which physical coordinates $ (x,y,z) $ are transformed into natural coordinates $ (r,s,t) $. The equivalence between these systems of coordinates in the local matrix calculation (see the Eq. ({\ref{Eq.7}})) is represented by: 

\begin{equation}
\begin{array}{lll}
\vec{k}_e & = & \displaystyle c \int \limits_{\Omega^e} \vec{B}^{T}(x,y,z) \vec{B}(x,y,z) dxdydz \\ & = & \displaystyle c \int \limits_{-1}^1 \int \limits_{-1}^1 \int \limits_{-1}^1 \vec{B}^{T}(r,s,t) \vec{B}(r,s,t) |\vec{J}(r,s,t)| drdsdt
\end{array}
\label{Eq.8}
\end{equation}
where $ \vec{J} $ is the Jacobian matrix and term $ |\vec{J}| $ is its determinant. After that, the Gauss quadrature is applied to the normalize domain (right hand side of Eq. ({\ref{Eq.8}})) as:

\begin{equation}
\vec{k}_e = \displaystyle c \sum_{i=1}^{M_i}  \sum_{j=1}^{M_j}  \sum_{k=1}^{M_k} \vec{B}^{T}(r_i,s_j,t_k) \vec{B}(r_i,s_j,t_k) |\vec{J}(r_i,s_j,t_k)| W_i W_j W_k
\label{Eq.9}
\end{equation}
where $ M_i, M_j, M_k, r_i, s_j, t_k $ and $ W_i, W_j, W_k $ are the number, points and weights of the integration points in each direction in a 3D domain. In our implementation, we use hexahedral FEs with eight nodes in the discretization of the domain, also known as 8-node brick elements which is a low-order tri-linear FE. The 8-node brick element was selected based on the idea that the majority of FE codes are based on the low-order elements and few works have focused on this type of elements {\cite{Knepley2013}}. The points and weights for an exact numerical integration of the 8-node brick element can be consulted in {\cite{Zienkiewicz2005}}.

\section{Numerical implementation}
\label{sec3} 
A typical way to construct the global sparse FE matrix is using serial computing. There are several ways to construct sparse matrices; however, the fastest method is to create the triplet\footnote{In the triplet sparse format the row and column indices that specify the position of NNZ (non-zero) entries must be stored. This means that each NNZ entry requires three numbers: the row and column position using integer numbers and the value of that entry using a double precision number. In this format, several NNZ entries with the same row/column indices can be repeated, which makes the difference with the coordinate (COO) sparse matrix format where entries with same row/column indices are summed.}  form first, and then, using a sparse function to compress the triplet sparse format, avoiding some memory challenges {\cite{Dabrowski2008}}. Despite the triplet format being simple to create, it is difficult to use in most sparse matrix algorithms; a compress sparse format is thus required {\cite{Davis2006}}. A common compress sparse format is the CSC (compressed sparse column), used by default in MATLAB {\cite{Gilbert1992}}. 

Algorithm {\ref{TradAssembly}} shows a typical implementation for the sparse matrix generation using the MATLAB syntax. This algorithm needs a subroutine that makes the numerical integration for computing local matrix $ \vec{k}_e $ and local load vector $ \vec{f}_e $. However, this work is focused only on generating matrix  $ \vec{K} $; load vector $ \vec{F} $ is not considered for parallel computation (this term is not considered because it is computationally much less demanding {\cite{Kruzel2013, Banas2014}}). The numerical integration subroutine changes with the PDE, element type and basis functions; it is a function of the nodal coordinates ($ \vec{C} $) and any nodal fields such as forces or boundary conditions ($ \vec{BCs} $) and material properties ($ \vec{MP} $) of element $ e $ {\cite{Cecka2011}}. After computing local matrix $ \vec{k}_e $, which is stored in array $ \vec{K}_{all} $, a mapping from local to global degrees of freedom (DOFs) is required. The mapping is made considering connectivity matrix ($ \vec{E} $), which gives the global node numbers that compose an element, and the number of DOFs per node (DOFxn). As a result of the mapping function, the positions in the global sparse matrix of each entry of $ \vec{k}_e $ is obtained, i.e. row ($ \vec{Ind}_{row} $) and column ($ \vec{Ind}_{col} $) indices, which are stored in arrays $ \vec{K}_{rows} $ and $ \vec{K}_{cols} $, respectively. Finally, the global sparse matrix is obtained in a single call of an assembly function (e.g. the MATLAB \emph{sparse} function {\cite{Davis2006, Dabrowski2008}}). 

\begin{algorithm}
\caption{Global sparse FE matrix construction in serial computing.}
\label{TradAssembly}
\begin{algorithmic}[1]
\State Initialize $ \vec{K}_{rows} $, $ \vec{K}_{cols} $ and $ \vec{K}_{vals} $ to zero;
	\For {each element $ e $ in $ nel $ }
		\State $ \vec{k}_e \longleftarrow NumericalIntegration(\vec{C,MP,BCs})$;
		\State $ \vec{K}_{vals}(e,:) = \vec{k}_e(:) $;
		\State $ (\vec{Ind}_{row}, \vec{Ind}_{col}) \longleftarrow Mapping(\vec{E}(e),DOFxn)$;
		\State $ \vec{K}_{rows}(e,:) = \vec{Ind}_{row}(:) $;
		\State $ \vec{K}_{cols}(e,:) = \vec{Ind}_{col}(:) $;
	\EndFor
\State $ \vec{K} = sparse(\vec{K}_{rows}, \vec{K}_{cols}, \vec{K}_{vals}) $;
\end{algorithmic}
\end{algorithm}

Several methods have been proposed to parallelize the global sparse matrix construction. In this work, instead of viewing the matrix generation subroutine as a single block, it is separated in two subroutines: the numerical integration subroutine and the assembly subroutine, which allows the analysis looking for parallelism opportunities. The numerical integration subroutine has several opportunities for parallelism {\cite{Banas2014}}:
\begin{enumerate}
	\item Computing all Gauss integration points in parallel;
	\item Parallelizing each $ \vec{k}_e $ component {\cite{Kruzel2013, Banas2014}};
	\item Finally, it is also possible compute $ \vec{k}_e $ in parallel for all elements {\cite{Knepley2013, Banas2014}}. 
\end{enumerate}
The loop over integration points (first approach) is a good strategy for higher-order FEs but it is harder to parallelize due to inherent data race\footnote{The race condition or race hazard is the behavior of a system in which the output is dependent on a sequence of events. It becomes a bug when events do not happen in the proper order. Race conditions can occur in computer software, especially multithreaded or distributed programs.}  in updating entries in $ \vec{k}_e $ with contributions from different integration points, and the degree of concurrency is generally lower  {\cite{Banas2014}}. The second approach is also helpful for higher-order FEs, but not convenient for low-order elements. The parallelization of the loop over low-order elements (third approach) is the most natural way because a thread can be charged entirely of one local matrix $ \vec{k}_e $ {\cite{Banas2014, Cecka2011}}. Hence, the third approach is selected in our implementation since the 8-node brick FE offers low degree of parallelism using the first two approaches. 

In the computation of $ \vec{k}_e $ there are no dependences and just the information from the specific 8-node brick element is necessary. As a consequence, numerical integration for different elements are also perfectly parallelizable {\cite{Kruzel2013}}. A thread is used here in to compute each local matrix $ \vec{k}_e $. Then, with a one-dimensional grid of thread blocks, all FEs in the mesh can be computed. 

In the eight-node brick element for scalar problems (1 DOF per node), each $ \vec{k}_e $ is a square dense matrix having 64 ($ 8 \times 8 $) components, but exploiting the matrix symmetry, only the lower or upper triangular matrix part, including its diagonal, should be considered. This leads to a significant saving in floating point operations and memory requirements, since approximately half the local matrix components are needed, specifically, 36 components instead of 64. Additionally, there are calculations that need to be made only once for the same FE type, such as the shape function and their derivatives in the physical coordinates. We use the same FE type to discretize the whole domain. 

Algorithm {\ref{CUke}} presents the CUDA kernel to compute the local stiffness matrix $ \vec{k}_e $. A thread is responsible of computing the 36 components of $ \vec{k}_e $ and store them in the GPU global memory as a vector in array $ \vec{K}_{vals} $. For each integration point, each thread must serially perform several operations, such as the Jacobian matrix ($ \vec{J} $) computation, its determinant ($ |\vec{J}| $) and its inverse ($ \vec{J}^{-1} $), in addition to the gradient matrix ($ \vec{B} $) and the local stiffness matrix ($ \vec{k}_e $). The shape functions ($ \vec{N} $) and their derivatives with respect to natural coordinates ($ d \vec{N} $) evaluated at each integration point are computed only once and they are stored in a fast GPU memory. The computation of matrix $ \vec{J} $ requires three dense matrix-vector multiplications, i.e. the multiplications of $ d \vec{N} $ ($ 3 \times 8 $) by nodal coordinates in Cartesian coordinates ($ 8 \times 3 $). The computation of $ |\vec{J}| $ is performed with the determinant formula of a $ 3 \times 3 $ matrix and $ \vec{J}^{-1} $ with the formula to invert matrices of the same size. Matrix $ \vec{B} $ is computed as the multiplication of $ \vec{J}^{-1} $ ($ 3 \times 3 $) by $ d \vec{N} $ ($ 3 \times 8 $). Finally, element matrix $ \vec{k}_e $ is computed as presented in Eq. ({\ref{Eq.9}}) with some modifications to calculate only the lower triangular matrix components. 

\begin{algorithm}
\caption{CUDA kernel for numerical integration in the 8-node brick element.}
\label{CUke}
\begin{algorithmic}[1]
\State \_\_global\_\_ void $ NumericalIntegration( \vec{k}_e, \vec{E}, \vec{C}, \vec{MP} )$
	\State $ tid = blockDim.x * blockIdx.x + threadIdx.x $; // Thread ID
	\State // Initialize local variables  
		\If {$ tid < nel $ }
			\For {$ i=0; i<8; i++ $}	\  \ // Loop over integration points
			\State read $ \vec{N}(i) $ and $ d\mathbf{N}(i) $;
			\State compute $ \vec{J} $, $ |\vec{J}| $ and $ \vec{J}^{-1} $;
			\State compute $ \vec{B} $;
			\State compute the lower triangular part of $ \vec{k}_{e} $ and store it in $ \vec{K}_{vals} $;
			\EndFor	
	\EndIf
\end{algorithmic}
\end{algorithm}

However, if the FE mesh is denser, the numerical integration may not be executed because this is a heavily memory-consuming operation and the GPU memory cannot allocate it entirely {\cite{Cecka2011}}. Thus, a simple domain partitioning based on the GPU memory availability is used, which is a modified implementation of that proposed by Komatitsch et al. {\cite{Komatitsch2009}}. Consequently, the numerical integration is developed in subgroups of FEs ($ GPU_g $), i.e. the numerical integration kernel is invocated to operate a FEs set. The elements subgroups are executed sequentially by the GPU, but the numerical integration over FEs within a group are executed in parallel. With this approach, the construction of the global FE matrix can be made from any mesh size, only limited by the CPU memory. Thus, if the available GPU memory is $ GPU_{ma} $ and the required memory by the numerical integration subroutine is $ GPU_{mr} $ (known through GPU query functions), the work done by the GPU is divided into groups based on the following formula:

\begin{equation}
GPU_g = ceil\left( \frac{GPU_{mr}}{GPU_{ma}} \right) 
\label{Eq.10}
\end{equation}

Each group of FEs is integrated numerically in the GPU, and the results allocated in the global GPU memory are downloaded to the CPU memory, leaving the GPU memory free for a new group of FEs. Every kernel is uploaded with only the data required by the FEs belonging the group that will be executed (only part of arrays $ \vec{C} $ and $ \vec{E} $), which allows a better use of memory. 

After the numerical integration over all elements, the assembly of all local matrices into the global matrix is executed. To do this, the assembly step can be viewed as a sparse matrix format conversion {\cite{Dziekonski2012a, Dabrowski2008}}. As stated before, the global position of each component of the $ \vec{K}_{all} $ matrix should be stored in arrays $ \vec{K}_{rows} $ and $ \vec{K}_{cols} $, which is an implicit sparse matrix format storage {\cite{Davis2006}}. The assembly step consists in converting the triplet format into the CSC format. This function requires sorting and reducing operations, which are not well suited for parallel computing due to they are intrinsic serial operations. In fact, the assembly cannot be considered an \emph{embarrassingly} parallel algorithm with no dependencies, as is considered the numerical integration is {\cite{Kruzel2013}}. The dependencies present in the assembly algorithm can cause race condition, which should be avoided to prevent wrong results. 

To remove dependencies present in the assembly algorithm there are some options, but the best for GPU programming is the use of coloring techniques or atomic operations {\cite{Dziekonski2012a, Cecka2011}}. Using the first option, sets of elements with different colors are obtained, with two elements having the same color if their local stiffness matrices do not coincide at any global matrix entry. Given this element partition, the numerical integration algorithm followed by immediate matrix assembly can be executed in parallel for elements within a group of the same color and proceeds color by color sequentially. However, this option requires pre-computing operations to divide the mesh into colors, where the problem to find the optimal number of colors is NP-hard. Although there are many heuristics for finding nearly minimal colorings {\cite{Kubale2004}}, this operation can take substantial execution time in meshes with a large number of FEs. On the other hand, atomic operations are undesirable because the penalty on the GPU execution time can be high {\cite{Cecka2011}}.

Dziekonski et al. {\cite{Dziekonski2012a}} use the sparse matrix conversion concept, which is executed in the GPU with the advantage that the assembly occurs in parallel using atomic operations. Using this approach, arrays $ \vec{K}_{rows} $ and $ \vec{K}_{cols} $ must be computed and stored in GPU, which substantially increases memory requirements. Although this approach is good on high-end GPUs provided with high-memory capacity to avoid memory transactions between processors, it is less interesting for many low-end video cards users. Additionally, in many current personal computers, the memory installed in the CPU motherboard is often four to eight times as large as than that installed on the graphics card {\cite{Komatitsch2009}}. Hence, a good option is to make the sparse matrix format conversion in the CPU.

\begin{figure}[hbtp]
\centering
\includegraphics[width=1\textwidth]{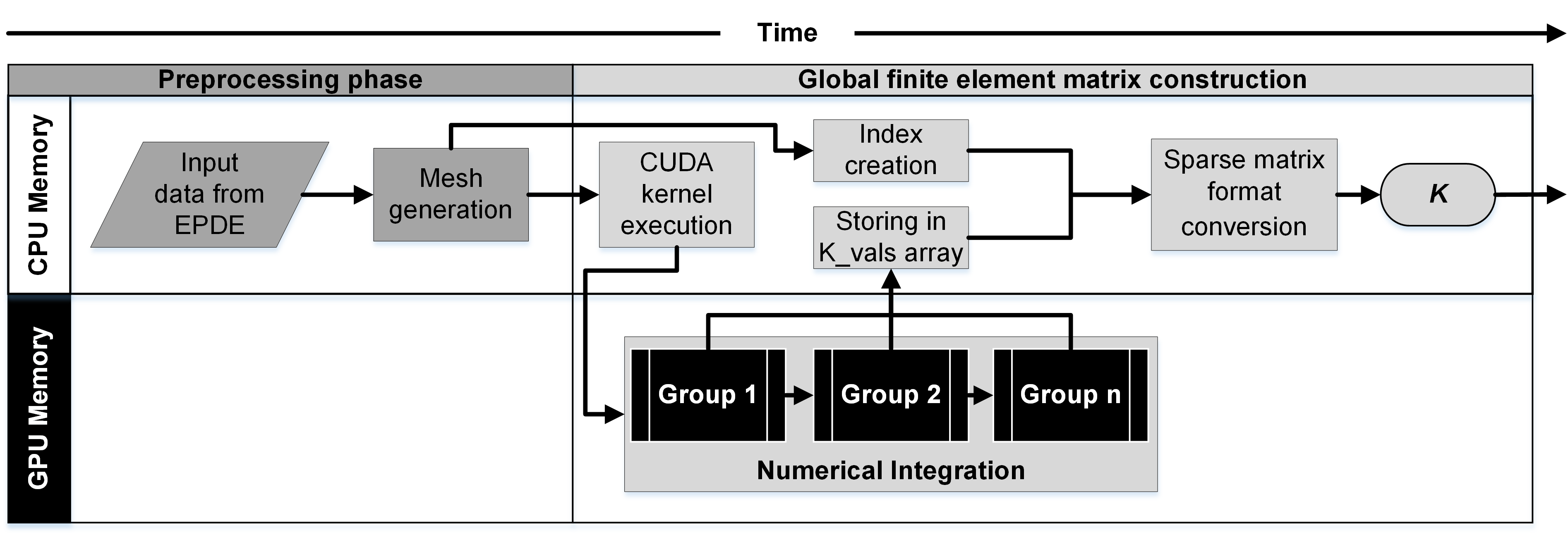}
\caption{Global sparse FE matrix construction.}
\label{MatGenSch}
\end{figure}

Therefore, a methodology for using GPU for intrinsic parallel and CPU for intrinsic serial algorithms  in the global FE matrix construction is proposed. The components of $ \vec{K}_{vals} $  are parallel computed in GPU, but indices $ \vec{K}_{rows} $ and $ \vec{K}_{cols} $ are   computed and stored in the CPU, since their computation is simple {\cite{Dabrowski2008}}.    After computing the global FE matrix in triplet format in CPU and GPU, as described above, the sparse matrix format conversion is executed in MATLAB using one of the following options: the native MATLAB function \emph{sparse}\footnote{MATLAB R2013b, \url{http://www.mathworks.com/help/matlab/ref/sparse.html}}, the \emph{sparse2} function from SuiteSparse\footnote{SuiteSparse 4.0.2, \url{http://www.cise.ufl.edu/research/sparse/SuiteSparse/}} and the \emph{sparse\_create} function from MILAMIN\footnote{Mutils 0.3, \url{http://www.milamin.org/}}. A graph showing the flow of our proposed implementation is presented in Fig. {\ref{MatGenSch}}. 

In Fig. {\ref{MatGenSch}}, there is a preprocessing phase that can be performed by any specialized software for complex domains. However, as our goal is evaluating the computational time involved in the matrix generation, we use cubic domains discretized by 8-node brick elements. This simple discretization is performed in MATLAB generating three arrays: connectivity matrix $ \vec{E} $, which indicates the nodes that compose an element; nodal matrix coordinates $ \vec{C} $, which specify the Cartesian coordinates of the nodes; and material properties matrix $ \vec{MP} $, which states the material properties of each element. As this subroutine is simple and takes a few seconds to generate large meshes, it is executed in the CPU. 

The second phase in Fig. {\ref{MatGenSch}} shows the necessary steps for building the global FE matrix, as described previously. In this step, the heterogeneous computing concept is used: while the GPU is executing the numerical integration kernel, the CPU is used for generating the indices of the sparse matrix in triplet format. Once the local matrices from all elements are computed, the assembly function can be executed, and the global matrix is obtained. After the global system is obtained, the boundary conditions are applied, and the solution of the linear system is found, but these steps are not shown in the figure because it is out of the scope of this paper. 

\section{Results}
\label{sec4}
The GPU results are obtained using a 2GB NVIDIA Quadro 4000 video card based on Fermi architecture. The CPU computations were performed in MATLAB R2013b using an Intel Xeon E5-2620 CPU at 2.0 GHz with six cores in a 64-bit Windows 7 operating system with 40GB of RAM. However, in the serial implementation only one CPU core is used. As stated before, the matrix is generated in two steps: first, conducting the numerical integration (NI) over all FEs, and second, assembling all local matrices into a global matrix. The results are presented separately in the following subsections. 

\subsection{Numerical integration}
\label{sec4.1}
The NI subroutine is programmed in CUDA C and is called from MATLAB workspace using the \emph{parallel.gpu.CUDAKernel} constructor, which creates a \emph{CUDAKernel} object. The properties of the CUDAKernel object such as the number of blocks in the grid and the number of threads in the block are also configured in MATLAB. After the creation and configuration of the kernel, it is evaluated using the \emph{feval} MATLAB function. 

\begin{figure}[hbtp]
\centering
\includegraphics[width=1\textwidth]{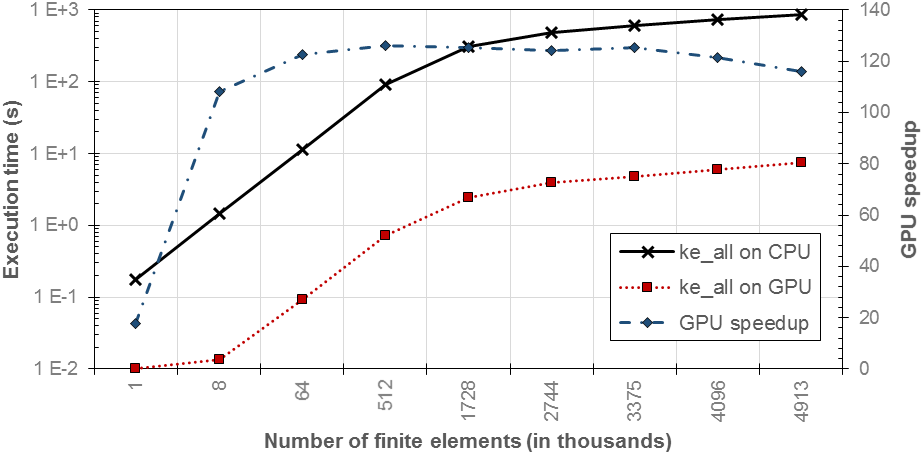}
\caption{Execution time in the numerical integration subroutine and GPU speedup as a function of the mesh size.}
\label{NIs_rst}
\end{figure}

The code executed in GPU is compared with its respective serial implementation counterpart executed in one CPU core. Fig. {\ref{NIs_rst}} shows the computational time when the NI routine is executed in CPU and GPU for different meshes. In the time spent by the GPU, the time spent in memory transactions is also taking into account: uploading data from CPU to GPU memory using the MATLAB function \emph{gpuArray}, and vice versa, downloading data from GPU to CPU memory using the function \emph{gather}. The gather function guarantees that all work is finished in the GPU, and it ensures accurate time measurement using the \emph{tic-toc} MATLAB functions. Fig. {\ref{NIs_rst}} also shows the speedup reached by the GPU implementation. The GPU speedup is calculated as the ratio between the serial ($ t_s $) and parallel ($ t_p $) execution time. The advantage of the GPU implementation is observed in the curves presented. 

The NI routine can compute meshes with almost five million of 3D FEs in a GPU with-low-memory capacity without using mesh partitioning. Although the NI routine is considered memory bounded {\cite{Cecka2011}}, the implemented code consumes little memory due to the low-order of the FE selected. The maximum time spent in the GPU for the largest mesh is 7.56 seconds, while for the same mesh in serial computing in the CPU, the NI routine spent 875.1 seconds, which means a GPU speedup of 115.8X. The GPU speedup dramatically increases its value using 8000 FEs. After that, the GPU speedup rate is slow; indeed, with more than four million FEs, the GPU speedup decreases. This fact is explained by the ratio between floating point operations (flops) and memory operations (memops). When the data volume transferred from one to another memory type (memops) is larger than the flops; the problem thus becomes memory bounded and the speedup is decreased. This is known as memory wall, in which there is a difference between the speed at which data are transferred to the processor and the rate at which instructions are executed {\cite{Asanovic2006}}. A common problem in transferring data between the CPU and the GPU memories is found in the PCI-express port, which limits the rate of moving data between these processors. For this reason, minimizing the transaction between the processor memories is recommended.

\subsection{Assembly}
\label{sec4.2}
The assembly step is performed in the CPU for two reasons: the assembly algorithm presents few parallel opportunities and requires high temporary memory. As the CPU cores are designed to maximize the execution speed of sequential programs (latency\footnote{Latency, in a basic definition, is the amount of time to complete a task, which mean that it is measured in units of time. While the CPU tries to minimize the latency, the GPU tries to maximize the throughput. Throughput is tasks completed per unit time measured in units as stuff per time, such as jobs completed per seconds. Therefore, the NI subroutine is executed in the GPU because we need to obtain many local matrices in the unit of time; and the assembly subroutine is executed in the CPU because the whole task is required to be finished in a short time. }) {\cite{Kirk&Hwu2012}}, this processor is well suited for serial tasks, as proposed here (assembly). The memory requirements of the assembly routine are shown in Table  {\ref{Trip2CSC}}. 

As indicated in Section {\ref{sec3}}, the assembly phase can be viewed as a sparse matrix conversion from triplet to CSC format. Table {\ref{Trip2CSC}} shows the percentage of compression in the non-zero elements (NNZ) going from one to another sparse format through different meshes. This table also shows the percentage of memory saving when the CSC is used instead of the triplet form. The percentage of compression and memory saving are not equal due to the matrix structure, which depends on the nodes and elements numeration; however, these percentages seem to be very similar. 

Having into account that the assembly subroutine initially requires a lot of memory for storing the three arrays  $ \vec{K}_{rows} $,  $ \vec{K}_{cols} $ and  $ \vec{K}_{vals} $, it is not feasible for a GPU implementation, particularly if GPUs with low-memory capacity are used for forming large matrices. 

\begin{table}[htb]
\begin{center}
\caption{Memory requirements in the triplet and CSC sparse formats for different mesh discretizations.}
\resizebox{\textwidth}{!} 
{
\begin{tabular}{l| c c c c c c c c c}
\hline 
\# of FEs & 1\,000 & 8\,000 & 64\,000 & 512\,000 & 1\,728\,000 & 2\,744\,000 & 4\,096\,000 & 5\,832\,000 & 8\,000\,000 \\ 
\hline 
NNZ in Triplet & 36\,000 & 288\,000 & 2\,304\,000 & 18\,432\,000 & 62\,208\,000 & 98\,784\,000 & 147\,456\,000 & 209\,952\,000 & 288\,000\,000 \\ 
NNZ in CSC & 15\,561 & 118\,121 & 920\,241 & 7\,264\,481 & 24\,408\,721 & 38\,710\,841 & 57\,728\,961 & 82\,135\,081 & 112\,601\,201 \\ 
NNZ compression & 56.8\% & 59.0\% & 60.1\% & 60.6\% & 60.8\% & 60.8\% & 60.9\% & 60.9\% & 60.9\% \\ 
\hline 
Memory used (triplet) & 0.58 MB & 4.61 MB & 36.9 MB & 294.9 MB & 995.3 MB & 1\,580.5 MB & 2\,359.3 MB & 3\,359.2 MB & 4\,608.0 MB \\ 
Memory used (CSC) & 0.26 MB & 1.96 MB & 15.3 MB & 120.5 MB & 404.7 MB & 641.8 MB & 957.0 MB & 1\,361.6 MB & 1\,866.6 MB \\ 
Memory saving & 54.9\% & 57.4\% & 58.6\% & 59.1\% & 59.3\% & 59.4\% & 59.4\% & 59.5\% & 59.5\% \\ 
\hline 
\end{tabular} 
}
\label{Trip2CSC}
\end{center}
\end{table}

\begin{figure}[hbtp]
\centering
\includegraphics[width=1\textwidth]{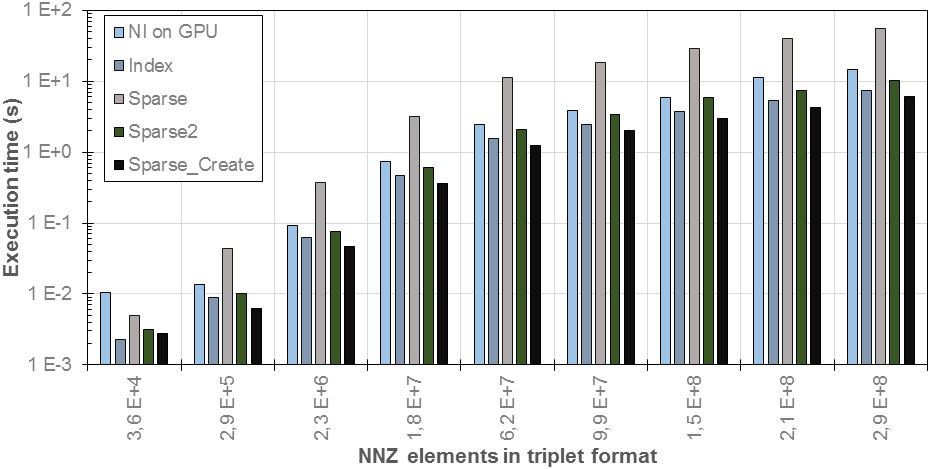}
\caption{Execution time spent by the assembly using three functions and different mesh sizes.}
\label{As_rst}
\end{figure}

Three different routines for changing the sparse matrix format are considered herein: the \emph{sparse} function from MATLAB (Sparse), the \emph{sparse2} function from SuiteSparse (Sparse2) and the \emph{sparse\_create} function from MILAMIN (SparseCreate). The execution time spent in the NI routine, executed in the GPU, and the assembly routine, executed in the CPU, are presented in Fig. {\ref{As_rst}}. Additionally, this figure shows the time spent for computing arrays $ \vec{K}_{rows} $ and $ \vec{K}_{cols} $ (Index) used by Sparse and Sparse2 routines. The NI routine spends more time than Sparse2 and SparseCreate routines for all mesh discretizations, even though the NI routine is executed in parallel in the GPU and the assembly routines are executed serially in the CPU. This occurs because the NI is an intensive routine that requires many operations between small matrices, such as multiplications, additions and inversions, whereas assembly only requires sorting and reducing operations {\cite{Dziekonski2012a}}.

The \emph{sparse} MATLAB function consumes significantly more time and memory than the other assembly functions, and for meshes larger than 200 thousand NNZs, it requires more time than the NI routine. The large memory usage explains the drawbacks of this function. In this routine, the row and column indices must be entered as double type arrays, which increases the memory requirements and performance overhead making this function the most inefficient assembly routine {\cite{Dabrowski2008}}. In turn, although the \emph{sparse2} function has an excellent performance, this routine requires computing and storing the index arrays, and their computations require almost the same time as the assembly function itself, meaning that the computational time spent by \emph{sparse2} can be duplicated if the time for computing $ \vec{K}_{rows} $ and $ \vec{K}_{cols} $ is considered. 

Regarding the \emph{sparse\_create} function for the sparse matrix conversion, in all discretization cases, the computational time invested in assembly is always the lowest. Indeed, this function employs less time than the index arrays efficiently computed as MILAMIN {\cite{Dabrowski2008}}. Moreover, the \emph{sparse\_create} function does not require the computation of the index arrays for the assembly step, saving memory and computations. However, this function is a specialized routine to form only matrices arising in FEM, which can represent a disadvantage when general sparse matrices are required.

\subsection{Partitioning the numerical integration}
\label{sec4.3}
Finally, as shown above, the numerical integration executed in parallel by using GPU and assembling the NNZ element based on \emph{sparse\_create} CPU function are the most efficient implementation.  However, without doing a mesh partitioning for the NI routine, the problem size that can be solved using low-memory capacity GPUs is limited. Hence, using an approach as the one presented in Section {\ref{sec3}} for partitioning the mesh based on available GPU memory, FE matrices can be constructed as large as CPU memory allows. In addition, as the CPU memory is usually greater than the GPU memory, larger matrices can be constructed than those created only in the GPU.

\begin{figure}[hbtp]
\centering
\includegraphics[width=1\textwidth]{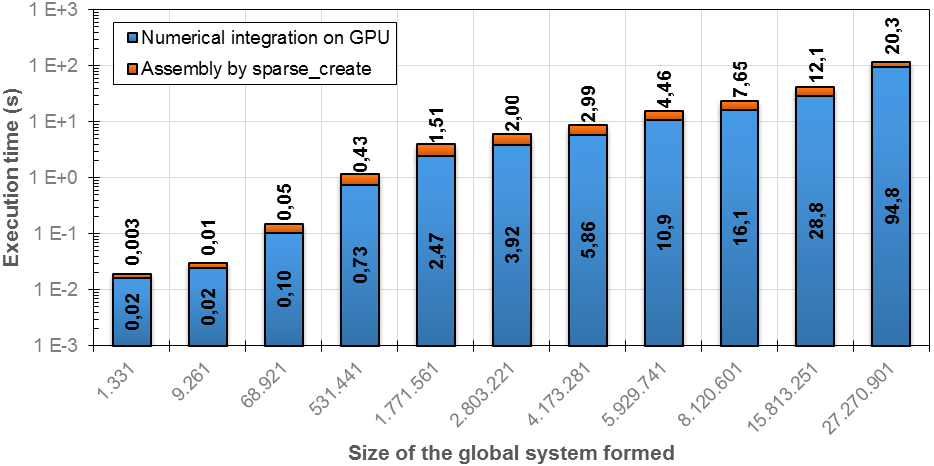}
\caption{Execution time spent generating global sparse FE matrices of several sizes.}
\label{PartNI_rst}
\end{figure}

Fig. {\ref{PartNI_rst}} shows the computational time for numerical integration and assembly routines in the construction of the global sparse matrix, considering different mesh sizes. Adding the time spent by these two routines, the time necessary to construct the global sparse matrix (MatGen) is obtained. Table {\ref{CTMatGen}} presents these results for different mesh sizes and the time percentage employed per function. The largest mesh size constructed has 27.27 million nodes and 27 million 3D elements. As stated before, the numerical integration routine always consumes more time than the assembly. Specifically, the assembly function consumes less than the 40\% of the total matrix generation time for all the mesh cases: 37.9\% is the greatest percentage spent for assembling 1\,728 million FEs (an intermediate mesh size) and the lowest percentage was presented assembling 1\,000 FEs (the lowest mesh size). These results indicate an appropriate approach to generate large FE matrices from unstructured meshes in systems in which the GPU memory is lower than the CPU memory. Additionally, this method shows the benefits of heterogeneous computing: using the CPU for intrinsic serial routines and the GPU for intrinsic parallel and intensive computations. 

\begin{table}[htb]
\begin{center}
\caption{Computational time spent in the sparse matrix construction.}
\resizebox{\textwidth}{!} 
{
\begin{tabular}{l |c c c c c c c c c c c}
\hline 
\# of FEs & 1\,000 & 8\,000 & 64\,000 & 512\,000 & 1\,728\,000 & 2\,744\,000 & 4\,096\,000 & 5\,832\,000 & 8\,000\,000 & 15\,625\,000 & 27\,000\,000 \\ 
\hline 
Matrix size & 1\,331 & 9\,261 & 68\,921 & 531\,441 & 1\,771\,561 & 2\,803\,221 & 4\,173\,281 & 5\,929\,741 & 8\,120\,601 & 15\,813\,251 & 27\,270\,901 \\  
MatGen time (s) & 0.02 & 0.03 & 0.15 & 1.17 & 3.98 & 5.92 & 8.85 & 15.3 & 23.8 & 40.9 & 115.0 \\  
NI time & 82.8\% & 80.9\% & 68.4\% & 62.9\% & 62.1\% & 66.2\% & 66.2\% & 70.9\% & 67.8\% & 70.5\% & 82.4\% \\ 
Assembly time & 17.2\% & 19.1\% & 31.6\% & 37.1\% & 37.9\% & 33.8\% & 33.8\% & 29.1\% & 32.2\% & 29.5\% & 17.6\% \\ 
\hline 
\end{tabular} 
}
\label{CTMatGen}
\end{center}
\end{table}

As a final comment, for reducing the time for matrix generation, one option is to form temporally small sparse matrices and next, to assemble all these matrices for forming the global matrix. Since the numerical integration for large meshes is divided into groups, which are executed sequentially in the GPU, it can can execute kernels asynchronously\footnote{Some function calls are asynchronous to facilitate concurrent execution, returning the control to the CPU thread before the GPU has completed a specific task. This means that the GPU is able to execute several kernels while the CPU executes other routines.}; an approach that can accelerate large sparse matrix generation in FEM is assembling all local matrices belonging to a FEs group while the GPU is executing the NI routine.

\section{Conclusions}
\label{sec5}  
The implementation proposed shows that large global sparse matrices arising in the numerical solution of EPDEs by FEM can be generated in heterogeneous computing with limited computational resources. In fact, with only 2GB of memory in the GPU, FEM matrices of 3D unstructured meshes with up to 27 million of FEs are generated. Moreover, a significant acceleration in the numerical integration executed in the GPU as compared to the serial CPU implementation is achieved. The maximum speedup reached by the GPU is 126x, and the mean value is 121x, which shows good performance of the developed code. Additionally, the assembly function (\emph{sparse\_create}) used in this work is an efficient routine that requires lower runtime and memory consumption compared with the other two functions (\emph{sparse} and \emph{sparse2}). In fact, \emph{sparse\_create} consumes less runtime than the numerical integration routine, even though the latter is implemented in parallel computing. The maximum time spent by the \emph{sparse\_create} function took only twenty seconds for assembling 27 million FEs, one-fifth of the time spent by NI routine.
 
Finally, although the proposed methodology is applied to the numerical solution of the Poisson equation by using FEM, this approach is general and can be used to construct not only the \emph{stiffness} matrix from static or stationary problems, but also the mass matrix from dynamics problems and for solving other types of PDEs by using the FEM.

%%%% Acknowledgments %%%%%%%%
\section*{Acknowledgments}
The first author thanks to Universidad Nacional de Colombia for the financial support to this work with the \emph{Estudiantes sobresalientes de posgrado} scholarship. MSG Tsuzuki was partially supported by CNPq (grant 310.663/2013-0).

\section*{References}

\bibliography{biblio}

\end{document}